\documentclass[authoryear,final,5p,times,twocolumn]{elsarticle}
\biboptions{sort&compress}

\usepackage{amssymb}
\usepackage{siunitx}
\usepackage[colorlinks=true,linkcolor=blue,urlcolor=blue,citecolor=blue,anchorcolor=blue]{hyperref}
\usepackage{physics} 
\usepackage{bm} 
\usepackage[capitalise]{cleveref}
\usepackage{tikz}
\usetikzlibrary{spy}
\usepackage{subcaption}
\usepackage{algorithmic}
\usepackage{xurl}
\usepackage{adjustbox}
\usepackage{pgfplots}
\usepackage{pdfpages}

\definecolor{TUDa-0d}{HTML}{535353}
\definecolor{TUDa-1a}{HTML}{5D85C3}
\definecolor{TUDa-2a}{HTML}{009CDA}
\definecolor{TUDa-3a}{HTML}{50B695}
\definecolor{TUDa-5a}{HTML}{DDDF48}
\definecolor{TUDa-8a}{HTML}{EE7A34}

\renewcommand{\mathbf}[1]{#1}

\journal{Engineering with Computers}

\begin{document}

\newcommand{\NurbsBasis}{\hat{N}}
\newcommand{\BSplineBasis}{\hat{B}}
\newcommand{\NurbsDegree}{p}
\newcommand{\NumberBasisFunctions}{n}
\newcommand{\NurbsWeight}{w}
\newcommand{\Jac}{\mathbf{J}_F}
\newcommand{\GeometryBasis}{G}

\newcommand\st{\mathrm{st}}
\newcommand\rt{\mathrm{rt}}
\newcommand\ag{\mathrm{ag}}
\newcommand\source{\mathrm{src}}

\newcommand{\MagVecPot}{\mathbf{A}}
\newcommand{\MagVecPotz}{A_z}
\newcommand{\MagVecPotzRt}{A_{z,\rt}}
\newcommand{\MagVecPotzSt}{A_{z,\st}}
\newcommand{\CurrentDensity}{\mathbf{J}}
\newcommand{\CurrentDensityTwoD}{J_{z,\source}}

\newcommand{\MagFluxDensity}{\mathbf{B}}
\newcommand{\Remanence}{\mathbf{B_\mathrm{r}}}
\newcommand{\RemanenceRed}{\mathbf{B^\bot_\mathrm{r}}}
\newcommand{\RotAngle}{\beta}
\newcommand{\MagnetAngle}{\alpha}
\newcommand{\CurrentAngle}{\varphi_0}


\newcommand{\ansatzFunction}{N} 
\newcommand{\ansatzFunctionIndex}{j}
\newcommand{\testSymbol}{v}
\newcommand{\testFunction}{N} 
\newcommand{\testFunctionIndex}{i}

\newcommand{\couplingSymbol}{G}
\newcommand{\couplingMatrix}{\mathbf{\couplingSymbol}}
\newcommand{\stiffness}{K}
\newcommand{\stiffnessMatrix}{\mathbf{\stiffness}}

\newcommand{\design}{x}
\newcommand{\adjoint}{z}
\newcommand{\state}{y}
\newcommand{\stateVector}{\mathbf{\state}}
\newcommand{\rhs}{b}
\newcommand{\rhsVector}{\mathbf{\rhs}}
\newcommand{\der}{\mathrm{d}}

\newcommand{\derSecond}{\mathrm{d}}

\newcommand{\Ckd}{C_\mathrm{kd}}
\newcommand{\Cle}{C_\mathrm{\ell e}}

\newcommand{\mortarSymbol}{\lambda}
\newcommand{\mortarVector}{\bm{\mortarSymbol}}

\newcommand{\opt}{\mathrm{opt}}
\newcommand{\fopt}{f_\opt}
\newcommand{\RotMat}{\mathbf{R}}
\newcommand{\ControlPoint}{C}
\newcommand{\OptiCtrlPoint}{\mathbf{\ControlPointSymbol}}
\newcommand{\Parameter}{P}
\newcommand{\OptiParameter}{\mathbf{\Parameter}}
\newcommand{\adjointSymbol}{\gamma}
\newcommand{\adjointVector}{\bm{\adjointSymbol}}
\newcommand{\Torque}{T}
\newcommand{\MeanTorque}{\overline{\Torque}}
\newcommand{\StdTorque}{\tilde{\Torque}}
\newcommand{\Amagnet}{A_\mathrm{Magnet}}
\newcommand{\Ttarget}{\Torque_\mathrm{Target}}

\newcommand{\OmegaRotor}{\Omega_\rt}
\newcommand{\OmegaStator}{\Omega_\st}
\newcommand{\GammaAirGap}{\Gamma_\ag}
\newcommand\IntegG{\mathrm{d}\Gamma}
\newcommand\Integ{\mathrm{d}\Omega}

\newcommand{\Iapp}{I_{\mathrm{app}}}
\newcommand{\nWind}{n_{\mathrm{wind}}}
\newcommand{\Acoil}{A_{\mathrm{coil}}}
\newcommand{\polePair}{p}

\definecolor{TUDa-2a}{HTML}{009CDA}
\definecolor{TUDa-2b}{HTML}{0083CC}
\definecolor{TUDa-3d}{HTML}{0071F3}
\definecolor{TUDa-3a}{HTML}{50B695}
\definecolor{TUDa-9b}{HTML}{E6001A}
\definecolor{TUDa-7b}{HTML}{F5A300}

\newcommand{\tbd}[1]{\textcolor{red}{TBD: #1}}
\newcommand{\mw}[1]{\textcolor{orange}{#1}}
	
\begin{frontmatter}
\title{Efficient Pareto-Front Generation for Electric Machines \\ using IGA and Second Order Derivatives}

\author[1]{Theodor Komann}
\author[2]{Michael Wiesheu}
\ead{michael.wiesheu@tu-darmstadt.de}
\author[1]{Stefan Ulbrich}
\author[2]{Sebastian Schöps}
\author[3]{Peter Gangl}

\affiliation[1]{organization={Department of Mathematics, Technical University of Darmstadt}, 
	city={Darmstadt},
	addressline={Dolivostraße~15}, 
	postcode={64297}, 
	country={Germany}}
\affiliation[2]{organization={Computational Electromagnetics Group, Technical University of Darmstadt},
	addressline={Schloßgartenstraße~8}, 
	city={Darmstadt},
    postcode={64289},
	country={Germany}}

 \affiliation[3]{organization={Johann Radon Institute for Computational and Applied
Mathematics, Austrian Academy of Sciences},
	addressline={Altenberger Straße~69}, 
	city={Linz},
    postcode={4040},
	country={Austria}}

\begin{abstract}
The multiobjective optimization of electric machines always involves a trade-off caused by various competing objectives such as performance and cost. A suitable design is usually determined by comparing variants from the Pareto front, which has been generated by a large number of simulation runs. This paper addresses the efficient generation of the Pareto front using a continuation method based on a homotopy method that exploits second-order derivative information to achieve superlinear convergence, enabling the fast generation of new Pareto-optimal points within only a few iterations. 
A key contribution is the derivation of formulas to compute the Hessian with respect to geometry parameters and shape, thus enabling direct modifications of the motor geometry in the context of Isogeometric Analysis. We apply our method to nonlinear 2D magnetostatic simulations of a permanent magnet synchronous motor and demonstrate its effectiveness by optimizing the cost, mean torque and torque ripple of the motor. 
Compared to a first-order optimization method, this approach reduces the number of iterations and function evaluations needed, making the pareto optimization fast and efficient.
\end{abstract}


\end{frontmatter}

\section{Introduction}
\label{sec:Introduction}
The increasing use of electric machines in transportation, particularly in the transition toward electric vehicles, highlights the importance of optimization in electric machine design. This shift represents a major change in the automotive industry, driven primarily by the need to address climate change and enhance environmental sustainability. In the last decade, advances in computer simulation have shifted the development of electric machines from physical to virtual prototypes, resulting in significant cost and material savings while promoting sustainability \citep{Boglietti_2009aa}. However, the current trend is to move beyond conventional numerical simulation methods and instead apply optimization strategies that can handle multiple design variables to improve performance.

In recent years, a variety of gradient-based approaches have been successfully developed for optimizing electric machines. Advances in free-form optimization have contributed in particular to the reduction of torque ripple \citep{Gangl_2015aa,Merkel_2021ab}, while parameter optimization has contributed to the reduction of magnet mass or maximizing torque \citep{Choi_2012aa,Kuci_2018aa,Lass_2017aa}. Furthermore, there is also work on the combination of the two optimization strategies \citep{Brun_2023,wiesheu_2024aa,wiesheu_2024ab}, to allow for a larger design space. Shape optimization has also been successfully applied in other engineering fields, such as fluid mechanics \citep{Soto2004,bletsos2021} and acoustics \citep{Schmidt2016,KAPELLOS2019368}. Although shape optimization alone typically dominates the computational cost, especially when many  degrees of freedom are involved, the simultaneous optimization of shape and parameters introduces a design space of  higher dimension. This makes the overall optimization more time-consuming than optimizing either the shape or the parameters alone. Optimizing sequentially, on the other hand, leads to a bilevel problem of the form $\min \min$, which requires deciding the order of parameter versus shape optimization. This sequential approach often leads to suboptimal solutions, as the full design space is not effectively explored in a stepwise way \citep{Sinha2018,dempe2002foundations}. By optimizing both simultaneously, we  overcome these limitations and achieve better designs. A table of existing design optimization methods applied for electric machines can be found in \cite[Table 1.3]{Brun_2024_thesis}

However, the challenges go beyond single-objective optimization. Conflicting criteria such as torque maximization versus cost minimization require multi-objective approaches.  From an industrial point of view, electric machines are  evaluated by several key performance indicators (KPIs) at once, such as average torque, torque ripple, losses and efficiency over multiple operating points, as well as mechanical and thermal performance. Improving one KPI typically worsens another, so finding the best compromise efficiently is challenging \citep{LIU2025115609,Sun_2021aa}. In this context, Pareto-optimal solution sets offer a systematic approach to industrial decision-making. They identify the possible compromises between competing KPIs, enabling engineers to select designs that best fulfil application-specific and economic requirements. This applies not only to electrical machines, but to many other engineering applications.  Despite this clear relevance, evolutionary algorithms (EAs), such as NSGA-II \citep{DebNSGA,Pereira2017}, have traditionally dominated multi-objective optimization of electrical machines \citep{Omar2022,Doi2019,Sato2022} due to their simplicity in application. However they are computationally expensive, i.e., need many evaluations of the underlying machine model. Gradient-based algorithms, on the other hand, can find optimized designs much faster, even in high-dimensional settings. To employ them in multi-objective settings, the vector-valued objective must first be scalarized, for example through the weighted-sum or epsilon-constraint formulations, which transform the problem into a series of scalar subproblems whose solutions approximate the Pareto front. In this context, several engineering disciplines have demonstrated the efficiency of adjoint-based approaches: in aerodynamic shape optimization, gradient-based adjoint methods achieve comparable Pareto front coverage to EAs at only a fraction of their computational cost \citep{Zingg2004}, while in electromagnetic design, adjoint-based homotopy methods further improve efficiency by  exploiting sensitivity information along the Pareto curve \citep{Gangl_2022ab,cesarano2024tracing}. Although these methods drastically reduce computation time, they typically provide only locally optimal solutions along the Pareto front rather than guaranteeing global optimality.

However the advantages of derivative-based multi-objective optimization are twofold. First, the use of the adjoint method enables computationally efficient derivative calculations, even for problems with many design variables. Unlike finite-difference or population-based approaches, the computational cost does not directly depend on the dimensionality. Second, gradients allow fast exploration of the design space, bypassing the ``curse of dimensionality'' that EAs suffer from. 
This is particularly important for modern electrical machine models, which have become increasingly complex in recent years due to the incorporation of hysteresis, mechanical and thermal effects \citep{Chen2020,Xiao2021,gangl2025multimaterialtopologyoptimizationelectric} in the FEM simulation. To trace Pareto-optimal solutions efficiently, gradient-based methods like predictor-corrector schemes are used. Initially introduced in \citep{Hillermeier_2001aa} and applied, e.g., in \citep{Martin2018}, they have also been adapted for shape optimization problems \citep{Gangl_2022ab,cesarano2024tracing}. However, a major drawback is their reliance on second-order derivatives, which are difficult and costly to implement in complex engineering applications. While automated differentiation  can provide these derivatives in theory, integrating AD into coupled multiphysics models which is common in electric machine design, often requires significant code modifications and computational power, limiting practical applications.
To our knowledge, existing work on the shape and parameter optimization of electrical machines uses first-order derivatives at best 
\citep{Brun_2023,wiesheu_2024aa,wiesheu_2024ab,Houta2023}. In contrast to recent work that demonstrates second-order shape optimization with standard FEM in an \textit{optimize-then-discretize} framework \citep{cesarano2024tracing,Gangl_2022ab}, our contribution focuses on a \textit{discretize-then-optimize} approach using IGA. We apply this to a state-of-the-art permanent magnet synchronous motor (PMSM) using the open source code GeoPDE \citep{Vazquez_2016aa}. We analytically derive second-order derivatives for both shape (control points) and parameter variables, extending the framework used in prior work \citep{wiesheu_2024aa,wiesheu_2024ab,Komann2024}. The use of second-order derivatives achieves local super-linear convergence of the design iterates towards a local minimizer, reducing the number of iterations required to reach a local optimal designs, even for highly nonlinear electric machine models. We exploit this property in multi-objective optimization to efficiently trace Pareto points via a predictor-corrector scheme.

The remainder of this paper is organized as follows:  \cref{sec:GeoRep} revisits the geometry representation based on spline functions and the physical model. \cref{Optimization} then formulates the optimization problem and derives first and second order derivatives with respect to shape and parameters. Numerical results for the second order optimization and Pareto front generation are presented in \cref{Results}. Finally, \cref{Conclusion} concludes the paper and discusses topics for future research.

\section{Modeling}
\label{sec:GeoRep}
\subsection{Geometry modeling}
When modeling complex geometries, choosing an appropriate representation method is important. Recently, spline-based techniques have become popular due to their ability to integrate simulation processes into Computer Aided Design (CAD) programs \citep{Nguyen2015}. Using a spline-based geometry technique comes with several advantages, particularly in the field of optimization, as demonstrated in \citep{Drzisga,Merkel_2021ab,WEEGER201926}. 
In regard to electric machines this is demonstrated in \citep{wiesheu_2024aa}, whose geometry model and methods also serve as foundation for this paper. 
We will now provide a brief overview of this concept which is called Isogeometric Analysis and  was introduced by Hughes and co workers in \citep{Hughes_2005aa}:

The univariate B-spline basis functions $N_{i}^{p}$ are defined by a knot vector
\[
    \Xi = \{\xi_1, \xi_2, \ldots, \xi_{n + p + 1}\},
\]
where each knot $\xi_i$ determines the support for the basis functions. The recursive construction of these functions is governed by the Cox-de Boor recursion formula \citep{DEBOOR1972}
\begin{equation}
 N_{i}^{p}(\xi) = \frac{\xi - \xi_i}{\xi_{i+p} - \xi_i} N_{i}^{p-1}(\xi) + \frac{\xi_{i+p+1} - \xi}{\xi_{i+p+1} - \xi_{i+1}} N_{i+1}^{p-1}(\xi),
 \label{eq:Cox-deBoor}
\end{equation}
with the recursion initialized as
\[
    N_{i}^{0}(\xi) = 
    \begin{cases}
    1 & \text{if } \xi_i \leq \xi < \xi_{i+1}\\
    0 & \text{otherwise}.
    \end{cases}
\]
The B-spline basis functions have several fundamental mathematical properties, including partition of unity, compact support, and non-negativity, which make them suitable for approximating geometries and shape optimization \citep{Piegl_1997aa,Cohen_2001aa}. 

\begin{figure*}[t]
    \centering
    \scalebox{0.85}{
        \input{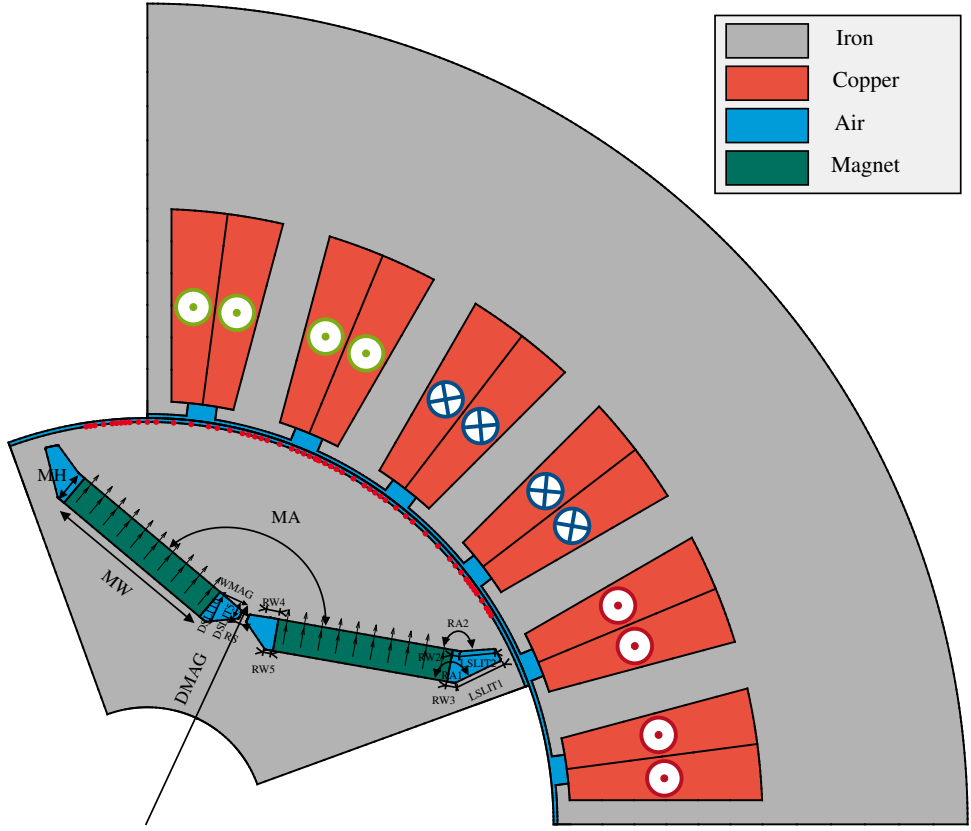}
    }
    \caption{Quarter model of the initial PMSM geometry including material definitions. Rotor parameters and control points that are used for optimization are highlighted.}
    \label{fig:geometry}
\end{figure*}
For higher-dimensional spline spaces, a tensor-product approach is applied. This technique involves multiplying individual basis functions from each parametric coordinate to form a set of multivariate basis functions. To map the parametric domain \( \hat{\Omega} \in (0,1)^d \) to the physical domain \( \Omega \subseteq \mathbb{R}^r \), where \( r \geq d \) we use Non-Uniform Rational Uniform B-Splines (NURBS) \( \hat{G}\), which generalize B-Splines to non-rational elements such as circles. Each geometric basis function \( \hat{G}_k \) is associated with a corresponding control point \( \mathbf{C}_k \). This results in the mapping \( \mathbf{\mathrm{F}}: \hat{\Omega} \rightarrow \Omega \)
\begin{equation}
    \mathbf{x}=\mathbf{\mathrm{F}}(\hat{\mathbf{x}}) = \sum\nolimits_{k} \hat{G}_k(\hat{\mathbf{x}}) C_k,
    \label{eq:MappingFromParametricToPhysical}
\end{equation}
The mapping  \eqref{eq:MappingFromParametricToPhysical} defines a diffeomorphism, which is required later to be able to calculate the derivatives with respect to the control points  \citep{Beirao-da-Veiga_2014aa}. For more details on the mathematical properties of isogeometric analysis  finite element discretization, see \citep{Bazilevs2006ISOGEOMETRICA}.

\subsection{Physical modeling}
\label{subsec:PhysicalModeling}
The state equation for our model is given by the magnetostatic low-frequency approximation of Maxwell's equations in 2D \citep{Jackson_1998aa}. This is the standard approach for the simulation of electric machines and described by 
\begin{equation}
    \begin{cases}
    \nabla \cdotp ( \nu \nabla \MagVecPotzRt) =\nu \nabla \cdotp \RemanenceRed & \mathrm{in}\ \OmegaRotor\\
    \nabla \cdotp ( \nu \nabla \MagVecPotzSt) =-\CurrentDensityTwoD & \mathrm{in}\ \OmegaStator,
    \end{cases} \label{eq:StrongRotorStator}
\end{equation}
as outlined before \citep{Bontinck_2018ac,Salon_1995aa}. We solve for the magnetic vector potential $A_z$ in rotor and stator domain including (nonlinear) material behavior given by the reluctivity $\nu$, which may depend on the magnitude of the magnetic flux density $B=\|\nabla A_z\|$. 
The right hand side is given by the $z$-component of the electric current density $J_{z,\source}$ for the stator and by the two dimensional remanence $\mathbf{B}_{r}$ of the permanent magnet  in the rotor \citep{Jannsen}. To model rotation, a mortaring approach is utilized, as described in \citep{De-Gersem_2004ad,Bontinck_2018ac}. 
For discretization we employ IGA, as outlined in \cref{sec:GeoRep}. We fix the degree $p=2$ and drop the superscript in $N^p_i$ for readability and follow the standard Ritz-Galerkin approach which results in
\begin{equation}
    \left(
    \begin{array}{ccc}
        \mathbf{K}_{\rt} & \mathbf{0} & -\mathbf{G}_{\rt} \\
        \mathbf{0} & \mathbf{K}_{\st} & \mathbf{G}_{\st} \mathbf{R}_{\alpha} \\
        -\mathbf{G}_{\rt}^{\top} & \mathbf{R}_{\alpha}^{\top} \mathbf{G}_{\st}^{\top} & \mathbf{0}
    \end{array}
    \right)
    \left(
    \begin{array}{c}
        \state_{\rt,\alpha} \\
        \state_{\st,\alpha} \\
        \lambda_{\alpha}
    \end{array}
    \right)
    =
    \left(
    \begin{array}{c}
        \mathbf{b}_{\rt,\alpha} \\
        \mathbf{b}_{\st,\alpha} \\
        \mathbf{0}
    \end{array}
    \right)
    \label{eq:CoupledPossionProblem}
\end{equation}
with the matrix entries
\begin{align}
K_{ij} & =\int _{\Omega } \nu(B) \nabla N_{i}(\mathbf{x}) \cdotp \nabla N_{j}(\mathbf{x}) \ \operatorname{d} \Omega, \label{eq:methodology:Krt} \\ 
b_{i,\rt} & =\int _{\Omega } \nu \nabla N_{i}(\mathbf{x}) \cdotp \mathbf{B}_{r}^{\bot } \ \operatorname{d}\Omega, \label{eq:methodology:brt} \\ 
b_{i,\st} & =\int _{\Omega } J_{z,\source} N_{i}(\mathbf{x})  \ \operatorname{d}\Omega. \label{eq:Methodology:bst}
\end{align} 
We denote the stiffness matrices as $\stiffnessMatrix_\rt$ and $\stiffnessMatrix_\st$ for rotor and stator, respectively. Note that $B_{r}$ and $J_{z,\source}$ vanish outside the magnets and windings, respectively, restricting the integration in \cref{eq:methodology:brt} and \cref{eq:Methodology:bst} to those subdomains. The solution vector $\state=(\state_{\rt,\alpha},
\state_{\st,\alpha},
\lambda_{\alpha})$ describes the state of the system at the rotation angle $\alpha$. 

The dependency on the control points becomes clear when performing the transformation to the reference domain $\hat{\Omega}$ using the transformation rule $\nabla N_{i} =\mathbf{J}_{F}^{-\top }\hat{\nabla }\hat{N}_{i}$ for functions in $H^1$ \citep{Monk_2003aa}. For example, for the stiffness matrix this results in
\begin{align}
    K_{ij} &=\int_{\hat{\Omega }} \nu(B) \left(\mathbf{J}_{F}^{-\top }\hat{\nabla }\hat{N}_{i}(\hat{x})\right) \cdotp \left(\mathbf{J}_{F}^{-\top }\hat{\nabla }\hat{N}_{j}(\hat{x})\right) |\mathbf{J}_{F}| \ \operatorname{d}\hat{\Omega }, \label{eq:methodology:Kparam} 
\end{align}
where $\mathbf{J}_{F}$ is the Jacobian matrix for the mapping \eqref{eq:MappingFromParametricToPhysical} and calculated with 
\begin{equation}
    \mathbf{J}_{F}= \begin{pmatrix}
    \sum\nolimits _{k}\left(\hat{\nabla }\hat{G}_{k}\right)^{\top }C_{kx}\\
    \sum\nolimits _{k}\left(\hat{\nabla }\hat{G}_{k}\right)^{\top }C_{ky}
    \end{pmatrix}.\label{eq:Optimization:JacobianMapping}
\end{equation}

Entries for the coupling matrices $\couplingMatrix_\rt$ and $\couplingMatrix_\st$ arise from harmonic (Fourier) functions used for mortaring on the interface of rotor and stator, given by
\begin{align*}
g_{i,2n-1} & =\int _{\Gamma _{\ag}} N_{i}\sin( n\theta )  \ \operatorname{d} \Gamma  & n\geq 1\\
g_{i,2n} & =\int _{\Gamma _{\ag}} N_{i}\cos( n\theta )  \ \operatorname{d} \Gamma  & n\geq 0
\end{align*} 
Efficient evaluations for different rotation angles are achieved by multiplying the stator coupling matrix with a rotation matrix $\RotMat_\alpha$, which contains block-diagonal sine and cosine entries. 
Overall, \eqref{eq:CoupledPossionProblem} is written as
\begin{equation}
\mathbf{e}(\state) = \stiffnessMatrix(\state)\state - \rhsVector = \mathbf{0}, \label{eq:Methodology:CoupledMatrixSystemred}
\end{equation}
where $\mathbf{e}$ is introduced as variable representing the state equation and the stiffness matrix depends of the solution vector, accounting for material nonlinearity. Here, the formula from \citep{Egger_2022ab} is used to calculate the torque via
\begin{equation}
T_{\alpha}(\state) = -l_z \state_{\st,\alpha}^{\top} \mathbf{G}_{\st} \mathbf{R}'_{\alpha} \lambda_{\alpha},
    \label{eq:Torque}
\end{equation}
where $l_z$ is the axial motor length.
To obtain the torque profile, we solve  \eqref{eq:Methodology:CoupledMatrixSystemred} at different rotor angles with a Newton-Raphson scheme and calculate the corresponding torque value with \eqref{eq:Torque}.

\section{Optimization}
\label{Optimization}
In this work, we follow a \textit{discretize-then-optimize} approach. This means that we choose a design variables within the finite-dimensional Euclidean space \citep{Wall_2008aa,Limkilde2021}. First, the geometry is built in  \(\mathbb{R}^2\) by parameters $\Parameter$, which define the general machine dimensions. Then, we select geometry control points $\ControlPoint$ as additional design variables for the optimization. The gradient with respect to these control points is calculated analytically using the adjoint formula \citep{Hinze_2009aa,troeltzsch2010optimal}. It has been shown that this approach leads to identical gradients compared to the \textit{optimize-then-discretize} approach in the context of shape optimization in IGA at least under some simplifying assumptions \citep{Fuseder_2015aa}. 

Building on this foundation, we now formalize a general framework for parameter and shape optimization within the IGA setting. The governing PDE is assumed to be of elliptic type, such as the magnetostatic problem introduced in \eqref{eq:CoupledPossionProblem}. This general formulation  will subsequently be then applied to the design of a PMSM.

\subsection{Combined parameter and shape optimization}
Following \citep{wiesheu_2024aa}, we  exploit the idea that changing a geometric parameter $\Parameter$ equals a simultaneous change of several control points $\ControlPoint$. 
If the derivatives of control point $\Ckd$ ($k$-th point in dimension $d$) are known, the respective derivative with respect to the $l$-th parameter $\Parameter_{l}$ can then be calculated with
\begin{equation}
  \frac{\der f}{\der \Parameter_{l}}
  \;=\;
  \sum_{k,d}
  \frac{\der f}{\der \Ckd}
  \;\frac{\der\Ckd}{\der \Parameter_{l}}.
  \label{eq:ParameterDerivative}
\end{equation}
This is necessary to optimize, for example, both the width of the magnet and the air gap simultaneously. To have a unified framework that can handle changes through both parameters and control points, each control point shift is interpreted as a parameter on its own.  

Assuming now, that also the second-order derivatives of an objective $f$ with respect to $\Ckd$ and $\Cle$ ($\ell$-th control point in dimension $e$) are known, by differentiating \eqref{eq:ParameterDerivative} again with respect to $\Parameter_{m}$, we obtain the second-order derivative
\begin{equation}
  \begin{aligned}
    \frac{\der ^{2} f}{\der \Parameter_{l}\,\der \Parameter_{m}}
    \;&=\;
    \frac{\der }{\der \Parameter_{m}}
    \Bigl(
      \sum_{k,d}\,\frac{\der f}{\der \Ckd}
      \;\frac{\der \Ckd}{\der \Parameter_{l}}
    \Bigr)
    \\[6pt]
    \;&=\;
    \sum_{k,d}\sum_{\ell,e}
    \frac{\der^{2} f}{\der \Ckd\,\der \Cle}
    \;\frac{\der \Cle}{\der \Parameter_{m}}
    \;\frac{\der \Ckd}{\der \Parameter_{l}}
    \;+\;
    \sum_{k,d}
    \frac{\der  f}{\der \Ckd}
    \;\frac{\der ^{2} \Ckd}{\der \Parameter_{l}\,\der \Parameter_{m}}. \label{eq:ParameterDerivativesSecond}
  \end{aligned}
\end{equation}
As in \eqref{eq:ParameterDerivative} the parameter derivatives can be seen as a postprocessing step. The second-order derivative of the control points with respect to the parameters \(\Parameter_{l}\) and \(\Parameter_{m}\) can be computed with second-order finite differences as in \citep{wiesheu_2024aa}. 
Note that \eqref{eq:ParameterDerivativesSecond} can also be applied in intermediate steps, e.g., when evaluating partial derivatives. This is beneficial when there are less design parameters compared to control points, which is usually the case. Therefore, we assume that \eqref{eq:ParameterDerivativesSecond} has been applied already for the partial derivatives in the following steps.

\subsection{Hessian matrix for the PDE-constrained problem}

Consider the optimization problem
\begin{equation}
\min f (\design, \state(\design)) \quad \text{s.t.} \quad \stiffness(\design,\state(\design))\state(\design) - \rhsVector(\design) = 0,
\label{eq:GeneralOptimizationproblem}
\end{equation}
where $f$ is the objective function we want to minimize.
Further, $\design$ denotes the optimization variables containing control points or  (geometric) parameters, $\state$ represents the state variables,  $\stiffness$ is the stiffness matrix and $\rhs$ is the right hand side. Since changing the design given by $\design$ changes the solution $\state$, $\state=\state(\design)$ is a function of $\design$. The variables $\stiffness$=$\stiffness(\design,\state(\design))$ and $\rhs=\rhs(\design)$ also depend on $\design$ and possibly $\state$. For ease of notation, these dependencies are not written out in the following.

To compute the second derivatives of the reduced functional $\bar{f}(\design) \coloneqq f(\design, \state(\design))$, where the state $\state(\design)$ is implicitly defined as the unique solution of the constraint in \eqref{eq:GeneralOptimizationproblem} efficiently, we start by introducing the Lagrangian
\[
L(\design,\state,\adjoint)= f(\design,\state) + \adjoint^\top e(\design,\state)
\]
under the assumption that all involved mappings are twice differentiable. The adjoint variable (Lagrange multiplier) is defined as \( \adjoint \). For the sake of concise notation in the subsequent derivation, we denote partial derivatives of the Lagrangian $L$ using subscripts. Specifically, we define the first partial derivative with respect to the design variables as $L_{\design} \coloneqq \frac{\partial L}{\partial \design}$ and the second partial derivative as $L_{\design\state} \coloneqq \frac{\partial^2 L}{\partial \design \partial \state}$. This convention applies analogously to the state and adjoint variables. To obtain the second derivative of the reduced functional $\bar{f}$ with respect to~$\design$, one has to proceed as follows \citep{Hinze_2009aa}:

\begin{enumerate}
    \itemsep0.5em 
    \item \textbf{Solve the State Equation:}
    We solve the nonlinear state equation
    \begin{equation}
        L_{\adjoint} = e(\design,\state) = \stiffness(\design,\state) \state - \rhs(\design) = \mathbf{0}
    \end{equation}
     iteratively with Newton-Raphson. This step is necessary for any nonlinear simulation. For the magnetostatic equations, the necessary derivatives are given in \ref{app:NewtonJacobian}.
    \item \textbf{Solve the Adjoint Equation:}
    The adjoint solution is obtained by solving
        \begin{align}
            L_{\state} &=0 & \Longleftrightarrow & & \left( \frac{\partial e(\design,\state)}{\partial \state} \right)^{\top} \adjoint &= -\left(\frac{\partial f(\design,\state)}{\partial \state} \right)^{\top}
        \end{align}
        for $\adjoint$. In the motor setting this corresponds to
        \begin{align*}
        \left( \left(\frac{\partial\stiffness}{\partial\state} \cdot\right)\state + \stiffness \right)^{\top} \adjoint &= -\left(\frac{\partial T_{\alpha}}{\partial \state} \right)^{\top}
        \end{align*}
        in order to calculate the torque derivatives for all rotation angles. The product is understood as:
        \begin{equation*}
            \left( \left(\frac{\partial \stiffness}{\partial \state} \cdot\right)\state \right)_{ij} = \sum_{k} \frac{\partial \stiffness_{ik}}{\partial \state_j} \, \state_k.
        \end{equation*}
    \item \textbf{Solve the Tangent Equation:} Next, we solve
        \begin{align}
         L_{\adjoint\state} \frac{\derSecond \state}{\derSecond \design} &= -L_{\adjoint\design}  &\Longleftrightarrow  & & \frac{\partial e(\design,\state)}{\partial \state} \frac{\derSecond \state}{\derSecond \design} = -\frac{\partial e(\design,\state)}{\partial \design}
        \end{align}
        for $\frac{\derSecond \state}{\derSecond \design}$. This is also called the direct sensitivity method that can be used as alternative to the adjoint method to obtain first order gradients. In our setting this corresponds to
        \begin{align*}
        \left( \left(\frac{\partial \stiffness}{\partial \state} \cdot\right)\state + \stiffness \right) \frac{\derSecond \state}{\derSecond \design} = -\frac{\partial \stiffness}{\partial\design} \state + \frac{\partial \rhs}{\partial\design}.
        \end{align*}

\item \textbf{Compute Hessian:}
After conducting the previous steps, we can calculate the Hessian with
\begin{align}
    \bar{f}_{\design\design} &= L_{\design\design} 
    + L_{\design\state} \frac{\derSecond\state}{\derSecond\design}
    + \left(\frac{\derSecond\state}{\derSecond\design}\right)^\top L_{\state\design} +
    \left(\frac{\derSecond\state}{\derSecond\design}\right)^\top L_{\state\state}\frac{\derSecond\state}{\derSecond\design} \\
    \Longleftrightarrow \quad \bar{f}_{\design\design} &= 
    \left( \frac{\partial^2 f(\design,\state)}{\partial \design^2} 
    + \adjoint^{\top} \frac{\partial^2 e(\design,\state)}{\partial \design^2} \right) \nonumber \\
    &\quad + \left( \frac{\partial^2 f(\design,\state)}{\partial \design\partial \state} 
    + \adjoint^{\top} \frac{\partial^2 e(\design,\state)}{\partial \design\partial \state} \right) 
    \frac{\derSecond \state}{\derSecond \design} \nonumber \\
    &\quad + \left(\frac{\derSecond \state}{\derSecond \design} \right)^\top
    \left( \frac{\partial^2 f(\design,\state)}{\partial \state\partial \design} 
    + \adjoint^{\top} \frac{\partial^2 e(\design,\state)}{\partial \state\partial \design} \right) 
     \nonumber \\
     &\quad + \left(\frac{\derSecond \state}{\derSecond \design} \right)^\top
    \left( \frac{\partial^2 f(\design,\state)}{\partial \state^2} + \adjoint^{\top} \frac{\partial^2 e(\design,\state)}{\partial \state^2} \right)
    \frac{\derSecond \state}{\derSecond \design} .
     \nonumber 
\end{align}

\end{enumerate}
The above procedure extends the first-order adjoint method to second-order. Overall, there are $n_f$ linear system solves necessary for the adjoint equation, and $n_\design$ linear solves  for the tangent equation. Yet, these evaluations are of moderate computational cost due to a single LU-decomposition of the state Jacobian and comparably small and sparse matrix systems ($n_\state<10000$). Currently, more expensive to evaluate is $L_{\design\design}$, as it involves the second derivatives of the stiffness matrix \eqref{eq:methodology:Krt} with respect to the design, which is defined by the control points. Recalling from \citep{Komann2024}, the first order derivatives with respect to the control points involve four integral contributions in the nonlinear case, see \ref{app:first-order-derivatives}. To obtain the second derivative
$L_{\design\design}$ analytically, the differentiation of those four terms results in a total of $20$ terms, which are given in \ref{app:second-order-derivatives:Lxx}. Although one can include the adjoint variable $\adjoint$ directly in the assembly process, $L_{\design\design}$ is still the most time-consuming part regarding computational cost. Also the derivative of $\rhs$ with respect to the control points needs to be calculated, where one obtains $6$ separate terms with the same argument. 
The derivatives $L_{\design\state}$ and $L_{\state\state}$ are comparably cheap, see \ref{app:second-order-derivatives:Lxy} and \ref{app:second-order-derivatives:Lyy}. For instance, computing 
\(\tfrac{\partial^2 T_{\alpha}}{\partial \state^2}\)  
yields a sparse matrix that is obtained by differentiating \eqref{eq:Torque} twice. 

One positive aspect of the heavy calculations for $L_{\design\design}, L_{\design\state}$ and $L_{\state\state}$ is that most of the steps involve only multiplication and addition, therefore the steps  are well suited for parallelization and GPU-computing \citep{Fatahalian2004}.

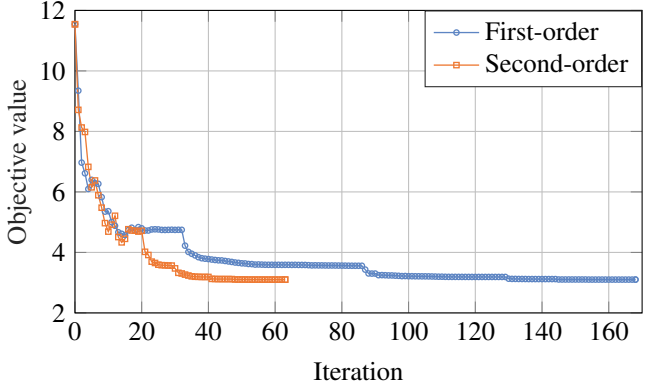
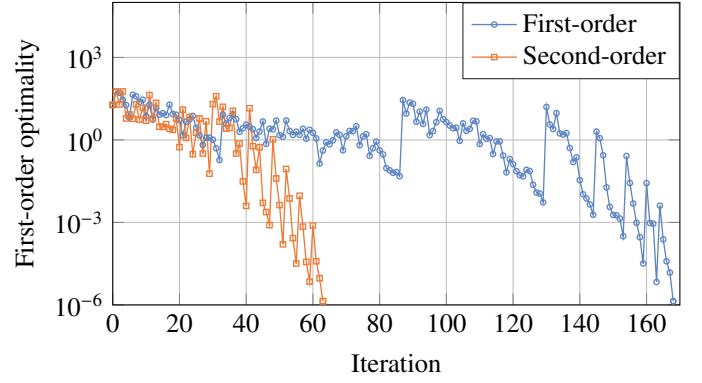
\begin{figure*}[t]
  \centering
    \begin{subfigure}[t]{0.5\textwidth}
        \centering
        \begin{tikzpicture}

\begin{axis}[%
    width=7.5cm,
    height=4cm,
    scale only axis,
    xmin=0,
    xmax=170,
    xlabel={Iteration},
    ylabel={Objective value},
    ymin=2,
    ymax=12,
    ylabel style={yshift=-0.5cm, font=\color{white!15!black}},
    xmajorgrids,
    ymajorgrids,
    legend style={at={(1,1)}, anchor=north east, legend cell align=left, align=left, draw=white!15!black}
]
\addplot [color=TUDa-1a, line width=0.5pt,mark=o, mark options={}, mark size=1pt]
table[x=iter,y=fval,col sep=comma] {data/FirstOrder.txt};
\addlegendentry{First-order}


\addplot [color=TUDa-8a, line width=0.5pt,mark=square, mark options={}, mark size=1pt]
table[x=iter,y=fval,col sep=comma] {data/SecondOrder.txt};
\addlegendentry{Second-order}

\end{axis}

\end{tikzpicture}
        \caption{Objective value versus iteration count.}
        \label{fig:FirstOrderIteration}
    \end{subfigure}%
    \begin{subfigure}[t]{0.5\textwidth}
        \centering
        \begin{tikzpicture}

\begin{axis}[%
    width=7.5cm,
    height=4cm,
    scale only axis,
    xmin=0,
    xmax=170,
    xlabel={Iteration},
    ylabel={First-order optimality},
    ymin=1e-6,
    ymax=1e5,
    ymode=log,
    ylabel style={yshift=-0.1cm},
    xmajorgrids,
    ymajorgrids,
    legend style={at={(1,1)}, anchor=north east, legend cell align=left, align=left, draw=white!15!black,nodes={scale=1.0}}
]
\addplot [color=TUDa-1a, line width=0.5pt,mark=o, mark options={}, mark size=1pt]
table[x=iter,y=firstorderopt,col sep=comma] {data/FirstOrder.txt};
\addlegendentry{First-order}


\addplot [color=TUDa-8a, line width=0.5pt,mark=square, mark options={}, mark size=1pt]
table[x=iter,y=firstorderopt,col sep=comma] {data/SecondOrder.txt};
\addlegendentry{Second-order}

\end{axis}

\end{tikzpicture}
          \caption{First‑order optimality versus iteration count.}
          \label{fig:FirstOrderOptimality}
    \end{subfigure}
  \caption{Comparison of first-order and second-order methods regarding objective value and first-order optimality over the iteration count. The initial objective decrease is similar, but after that the second-order method needs significantly less iterations to converge.}
  \label{fig:OptimizationFirstOrder}
\end{figure*}

\subsection{Application to a Permanent Magnet Synchronous Motor}
When optimizing PMSMs, there are several competing design criteria. Commonly used objectives for optimization are material cost, motor performance given by the average torque, and the quality of the torque profile (measured by torque ripple) \citep{Sato2022,COMSOL2021motoropt}. 
Here, we consider three objective functions to be minimized:
\begin{enumerate}
    \itemsep0.5em
    \item The rotor material cost, denoted by $M$.
    \item The torque standard deviation denoted by $\widehat{T}$, representing torque ripple.
    \item A smoothing (regularization) term $S$, that penalizes abrupt changes in the rotor surface, facilitating manufacturability.
\end{enumerate}
They are combined via a weighted sum approach with positive weights $m_{1}, m_{2}, m_{3}$. To maintain a certain performance, a lower bound for the average torque is prescribed. The resulting optimization problem then reads
\begin{equation}
  \label{eq:FullOptimizationProblem}
  \begin{aligned}
  \min_{\design} \;F(\design)
  \;=\;&
       m_{1}\,M(\design)
       \;+\;
       m_{2}\,\widehat{T}(\design)
       \;+\;
       m_{3}\,S(\design),
  \\[6pt]
  \text{subject to}\quad
  & \stiffness\left(\design,\state(\design)\right)\,\state
    \;-\;\rhs(\design)
    \;=\;\mathbf{0}
    \quad
    \text{(\ref{eq:CoupledPossionProblem})},
  \\[4pt]
  & \overline{T}
    \;\ge\;
    T_{\mathrm{Target}}
    \quad
    \text{(average torque constraint)},
  \\[4pt]
  & \mathbf{g}\bigl(\design\bigr)
    \;\le\;
    \mathbf{0}
    \quad
    \quad\text{(geometric feasibility)}.
  \end{aligned}
\end{equation}
where $\overline{T}$ is the average torque. The Term  $S(\design)$ is defined by \begin{equation*}
  \label{eq:S-Definition}
  S
  \;=\;
  \sum_{i}\,
    \frac{\bigl(\Delta C_{i+1} \;-\;\Delta C_{i}\bigr)^2}
         {\theta_{i+1} \;-\;\theta_{i}},
\end{equation*}
where $\theta_{i}$ and $\theta_{i+1}$ are consecutive angular positions of control points at the rotor surface, and $\Delta C_{i}$ measures the radial offset of the $i$-th point. The mean torque is defined by
\begin{equation}
  \label{eq:MeanTorqueDefinition}
  \overline{T}(\design)
  \;=\;
  \frac{1}{N}
  \sum_{\alpha=1}^{N}
      T_{\alpha}(\design).
\end{equation}
If the parameters in $\design$ are denoted by $P_{l}$, then by linearity,
\begin{equation}
  \frac{\mathrm{d}\,\overline{T}}
       {\mathrm{d}P_{l}}
  \;=\;
  \frac{1}{N}
  \sum_{\alpha=1}^{N}\,
      \frac{\mathrm{d}\,T_{\alpha}}
           {\mathrm{d}P_{l}}.
           \label{eq:MeanTorqueDerivative}
\end{equation}
We define the torque ripple by the standard deviation
\begin{equation}
  \widehat{T}
  \;=\;
  \sqrt{\frac{1}{N}
   \sum_{\alpha=1}^{N}
      \bigl(
         T_{\alpha} - \overline{T}
      \bigr)^2}.
\end{equation}
Its derivative with respect to $\Parameter_{l}$ follows from a standard chain-rule approach \citep{Li_2015ac}:
\begin{equation}
  \frac{\mathrm{d}\,\widehat{T}}
       {\mathrm{d}P_{l}}
  \;=\;
  \frac{1}{\widehat{T}}
  \biggl(
    \frac{1}{N}
    \,\sum_{\alpha=1}^{N}\,
       T_{\alpha}
       \;\frac{\mathrm{d}\,T_{\alpha}}{\mathrm{d}P_{l}}
    \;-\;
    \overline{T}
    \;\frac{\mathrm{d}\,\overline{T}}{\mathrm{d}P_{l}}
  \biggr).
  \label{eq:StdTorqueDerivative}
\end{equation}
More details about the optimization problem are found in \citep{wiesheu_2024aa}, here we focus on the details for higher order optimization. To apply a Newton method we now have to compute the second order derivatives of \eqref{eq:MeanTorqueDerivative} and \eqref{eq:StdTorqueDerivative}. Again due to linearity, differentiation of \eqref{eq:MeanTorqueDerivative} yields directly
\begin{equation}
  \label{eq:MeanTorqueParamSecondOrder}
  \frac{\mathrm{d}^{2}\,\overline{T}}
       {\mathrm{d}P_{l}\,\mathrm{d}P_{m}}
  \;=\;
  \frac{1}{N}
  \;\sum_{\alpha=1}^{N}\;
    \frac{\mathrm{d}^{2}\,T_{\alpha}}
         {\mathrm{d}P_{l}\,\mathrm{d}P_{m}}.
\end{equation}
The second derivative of $\widehat{T}$ is obtained by differentiating \eqref{eq:StdTorqueDerivative} and results in
\begin{align}
  \label{eq:SecondOrderStdTorque}
  \frac{\mathrm{d}^{2}\,\widehat{T}}
       {\mathrm{d}P_{l}\,\mathrm{d}P_{m}}
  &=\;
  -\,\widehat{T}^{-1}
  \biggl(\,
    \frac{\mathrm{d}\,\widehat{T}}{\mathrm{d}P_{m}}
  \biggr)
  \biggl(\,
    \frac{\mathrm{d}\,\widehat{T}}{\mathrm{d}P_{l}}
  \biggr)
  \nonumber\\[6pt]
  &\quad
  +\;\frac{1}{\widehat{T}}
  \biggl(
    \frac{1}{N}
    \Bigl[\,
      \sum_{\alpha=1}^{N}\;
         \frac{\mathrm{d}T_{\alpha}}{\mathrm{d}P_{m}}\,
         \frac{\mathrm{d}T_{\alpha}}{\mathrm{d}P_{l}}
      \;+\;
         T_{\alpha}\,
         \frac{\mathrm{d}^{2}T_{\alpha}}
              {\mathrm{d}P_{l}\,\mathrm{d}P_{m}}
    \Bigr]
    \biggr.
  \nonumber\\[6pt]
  &\quad
  \biggl.
    -\,\frac{\mathrm{d}\,\overline{T}}{\mathrm{d}P_{l}}\,
       \frac{\mathrm{d}\,\overline{T}}{\mathrm{d}P_{m}}
    \;-\;
    \overline{T}\,
       \frac{\mathrm{d}^{2}\,\overline{T}}
            {\mathrm{d}P_{l}\,\mathrm{d}P_{m}}
  \biggr).
\end{align}
We have now formulated a multi-objective optimization problem for electric machines and derived  necessary first- and second-order derivatives of the objectives with respect to the design parameters. In the next section, we compare the performance of a first-order gradient-based method with a second-order approach on the same problem to investigate their performance.

\section{Results}
\label{Results}
The presented optimization procedure has been implemented  in \textit{MATLAB}\textsuperscript{\textregistered}~R2024b relying on the \textit{GeoPDEs} toolbox \citep{Vazquez_2016aa} for the basic IGA functionality. As an application example, the PMSM shown in Figure \ref{fig:geometry} has been constructed using the NURBS toolbox \citep{NURBS_2021aa}. Also highlighted are the rotor parameters and control points of the rotor surface, that are allowed to change during optimization. The weights in the objective function \eqref{eq:FullOptimizationProblem} are set to $$m_1 = 1,\quad  m_2 = 100,\quad m_3 = 1000,$$
through all optimization steps. As the cost and the three terms of the objective differ by several orders of magnitude, the weights are used to bring all objectives to a similar range, thus ensuring that none of them dominates the optimization. All simulations are run on a standard laptop with (16-core Intel\textsuperscript{\textregistered} Core\textsuperscript{TM} Ultra 7@4.8GHz, 32GB RAM). There are 16 geometric parameters for optimization, 29 design variables that define the control point offsets and one variable for the electric phase, so overall we have $n_f=46$. The discretization of \eqref{eq:CoupledPossionProblem} is set to $n_y=6787$ degrees of freedom. More details on the discretization and choice of optimization variables can be found in \citep{wiesheu_2024aa}, as we follow a very similar setup. 

\begin{figure*}[t]
  \centering
    \begin{subfigure}[t]{0.5\textwidth}
        \centering
        \begin{tikzpicture}

\begin{axis}[%
    width=7.5cm,
    height=4cm,
    scale only axis,
    xmin=0,
    xmax=170,
    xlabel={Iteration},
    ylabel={First-order optimality},
    ymin=1e-6,
    ymax=1e5,
    ymode=log,
    ylabel style={yshift=-0.1cm},
    xmajorgrids,
    ymajorgrids,
    legend style={at={(1,1)}, anchor=north east, legend cell align=left, align=left, draw=white!15!black,nodes={scale=1.0}}
]
\addplot [color=TUDa-1a, line width=0.5pt,mark=o, mark options={}, mark size=1pt]
table[x=iter,y=firstorderopt,col sep=comma] {data/FirstOrder.txt};
\addlegendentry{First-order}

\addplot [color=TUDa-8a, line width=0.5pt,mark=square, mark options={}, mark size=1pt]
table[x=iter,y=firstorderopt,col sep=comma] {data/SecondOrder2_80.txt};
\addlegendentry{Second-order}

\end{axis}

\end{tikzpicture}
        \caption{Objective value versus iteration count.}
        \label{fig:MixedIteration}
    \end{subfigure}%
    ~
    \begin{subfigure}[t]{0.5\textwidth}
        \centering
        \begin{tikzpicture}

\begin{axis}[%
    width=7.5cm,
    height=4cm,
    scale only axis,
    xmin=70,
    xmax=170,
    xlabel={Iteration},
    ylabel={Objective value},
    ymin=3,
    ymax=3.6,
    ylabel style={yshift=-0.3cm, font=\color{white!15!black}},
    xmajorgrids,
    ymajorgrids,
    legend style={at={(1,1)}, anchor=north east, legend cell align=left, align=left, draw=white!15!black}
]
\addplot [color=TUDa-1a, line width=0.5pt,mark=o, mark options={}, mark size=1pt]
table[x=iter,y=fval,col sep=comma] {data/FirstOrder.txt};
\addlegendentry{First-order}

\addplot [color=TUDa-8a, line width=0.5pt,mark=square, mark options={}, mark size=1pt]
table[x=iter,y=fval,col sep=comma] {data/SecondOrder2_80.txt};
\addlegendentry{Second-order}

\end{axis}

\end{tikzpicture}
        \caption{Objective function evaluations versus iteration count.}
        \label{fig:MixedOptimality}
    \end{subfigure}%
  \caption{Comparison of the first-order method to a mixed approach, where the second-order continues the optimization. Since the optimum is already close, the second-order method converges rapidly.}
  \label{fig:MixedOptimization}
\end{figure*}
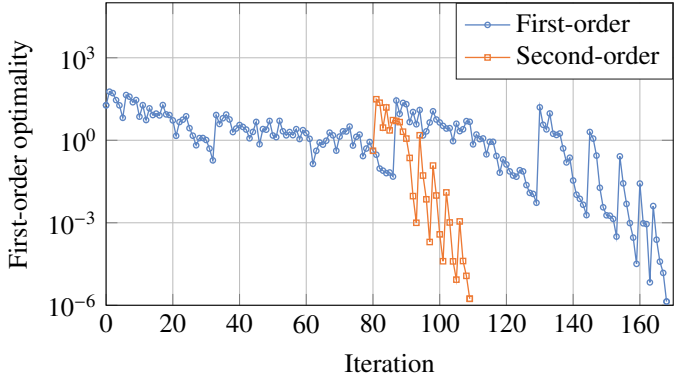
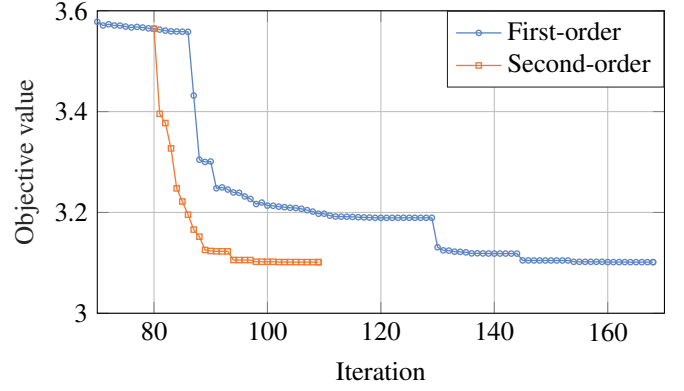

\subsection{Comparison of first- and second-order optimization}
\begin{figure}[h]
    \centering
    \includegraphics[width=0.75\linewidth]{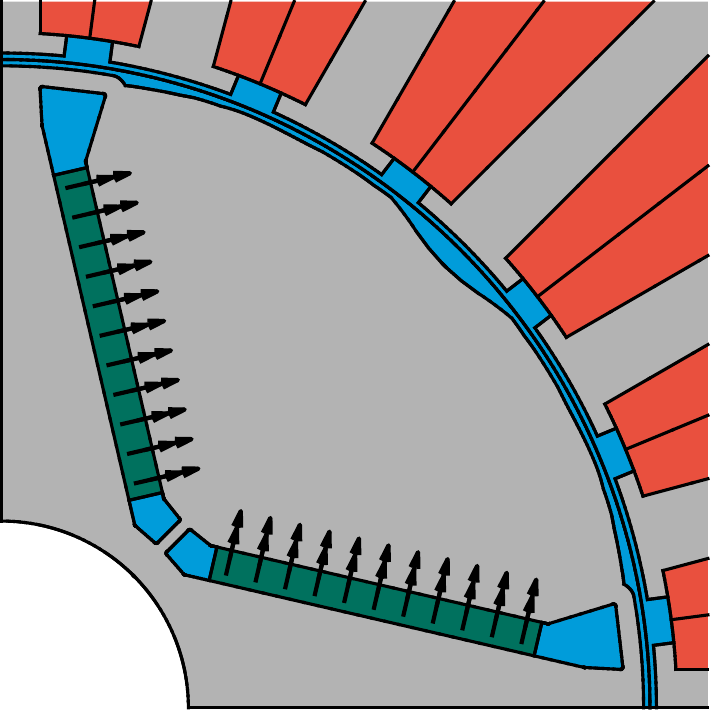}
    \caption{Optimized rotor geometry $T_\mathrm{Target}=\SI{0.575}{Nm}$.}
    \label{fig:OptimizationGeometry}
\end{figure}

We use the interior point solver \texttt{fmincon} from the optimization toolbox and provide the optimization objective, constraints, and their derivatives as defined before. For a first evaluation, we compare the optimization process for three different settings:
\begin{itemize}
  \item \textbf{First-order method:}
        We provide \texttt{fmincon} only gradients. Methodologically the same as in \citep{wiesheu_2024aa}.
  \item \textbf{Second-order method:} We additionally provide the exact Hessian matrix for the objectives and constraints. The barrier update is performed using the \textit{monotone} setting.
  \item \textbf{Mixed-order method:} We start using the first-order method for a fixed number of iterations. After that, we switch to the second-order method.
\end{itemize}
The target torque for the optimization is set to $T_\mathrm{Target} = \SI{0.575}{Nm}$ for the quarter model with a length of $l_z=\SI{3.5}{cm}$. The cost of permanent magnets and iron are set to dimensionless values of $50$ and $2$ per kg, respectively.
Both methods are provided with the same initial design from  \cref{fig:geometry} and identical termination settings, i.e., a first order optimality tolerance of $10^{-6}$. The resulting geometry for the optimization runs is shown in \cref{fig:OptimizationGeometry}. Both the parameters and the rotor surface are adjusted in order to reduce the magnet mass and the torque ripple during rotation.

A visualization of how the objective function evolves during the optimization process is given in \cref{fig:FirstOrderIteration}. As expected, the second-order method converges in less iterations compared to the first order. It reduces the number of iterations needed almost by two thirds (63 vs 168). The convergence behavior described by the first-order optimality is shown in \cref{fig:FirstOrderOptimality}.
Here, one can see the advantage of having second-order derivatives: Once the solver approaches the optimum, super-linear convergence is reached. Especially after iteration 40, whenever the barrier-term for the constraints is adjusted, the first-order optimality is drastically improved. 

\begin{table}[h]
    \caption{Comparison of iterations and computational runtime for the different optimization methods.}
    \label{tab:OptimizationComparison}
    \centering
    \renewcommand{\arraystretch}{1.15}
    \begin{tabular}{|l||c|c|c|c|}
    \hline
    \textbf{Method} & \textbf{Iterations} & \textbf{Evaluations} & \textbf{Runtime} \\ \hline\hline
    \textbf{First-order} & 168 & 244  & \SI{53}{min} \\ \hline
    \textbf{Second-order} & 63  & 101  & \SI{1}{h} \SI{39}{min} \\ \hline
    \textbf{Mixed-order}
           & 80+29 & 139+44  & \SI{1}{h} \SI{15}{min} \\ \hline
    \end{tabular}
\end{table}

When it comes to computational time, however, the first-order method still finishes earlier, as each Newton step is considerably more expensive: One iteration takes  approximately \SI{19}{s}, while one iteration for the second-order method needs about \SI{94}{s}.
In that time the gradient method has already completed multiple steps. Consequently, the total runtime of the second-order run (\SI{1}{h}\,\SI{39}{min}) surpasses the first-order runtime (\SI{53}{min}). However, it is expected that an efficiency-focused implementation for the second-order derivatives would resolve this issue. In addition, faster Hessian computations by using a Hessian vector product via dual adjoint solves could reduce the computation time of second order derivatives, as explicit assembly of the entire Hessian is not required \citep{enouen2022}. The computational times are summarized again in \cref{tab:OptimizationComparison}.

A natural idea to balance the low per‑iteration cost of the first-order approach with the local superlinear convergence of the second-order approach is to employ a switching strategy. In other words, one first uses the gradient‑only approach initially to reduce the objective quickly. Once the design is sufficiently progressed, one switches to a second‑order method to exploit its superlinear convergence near the optimum. This results in the mixed-order method. For instance, we performed 80 first-order iterations and then switched to second-order for 29 further iterations, leading to a total of 109 iterations. The resulting design is again the same but the computation time (\SI{1}{h} \SI{15}{min}) is reduced compared to the second-order method. This process is visualized in \cref{fig:MixedOptimization}.
After some initial iterations, the first-order optimality is rapidly improved and the solver converges in $29$ iterations instead of $88$ remaining iterations, see \cref{fig:MixedIteration}. Also the objective value is very quickly reduced as highlighted in \cref{fig:MixedOptimality}. For demonstration purposes, we selected 80 first-order iterations and then switched once the gradient solver had substantially reduced the objective. An automated strategy could track certain indicators (e.g. the relative change in the objective, the condition number of the approximate Hessian, or the KKT norm) to determine the best switching point. Even without these possible refinements, our basic experiment already demonstrates the potential of the mixed approach. Although still not faster than the first-order method, the point we want to make here is the advantageous convergence properties near an optimum, which will be exploited in the following.

\begin{figure*}[h]
    \centering
    \begin{subfigure}[t]{0.5\textwidth}
        \centering
        \begin{tikzpicture}

\begin{axis}[%
    width=7.5cm,
    height=4cm,
    scale only axis,
    xmin=0.45,
    xmax=0.7,
    xlabel={Torque (Nm)},
    ylabel={Cost},
    xtick={0.45,0.5,...,0.7},
    xlabel style={yshift=-0.cm},
    ylabel style={yshift=-0.35cm},
    xmajorgrids,
    ymajorgrids,
    legend style={at={(1,0)}, anchor=south east, legend cell align=left, align=left, draw=white!15!black}
]

\addplot [color=TUDa-1a, line width=1pt,mark=o, mark options={}, mark size=2.0pt]
table[x=torque,y=cost,col sep=comma] {data/Pareto1stInit.txt};
\addlegendentry{First-order pareto}

\addplot [color=TUDa-8a, line width=1pt,mark=square, mark options={}, mark size=2.0pt]
table[x=torque,y=cost,col sep=comma] {data/Pareto2ndInit.txt};
\addlegendentry{Second-order pareto}

\end{axis}

\end{tikzpicture}
        \caption{Comparison of Pareto fronts for first-order vs.\ second-order solutions.}
        \label{fig:pareto:cost}
    \end{subfigure}%
    ~
    \begin{subfigure}[t]{0.5\textwidth}
        \centering
        \begin{tikzpicture}

\begin{axis}[%
    width=7.5cm,
    height=4cm,
    scale only axis,
    xmin=0.45,
    xmax=0.7,
    xlabel={Torque (Nm)},
    ylabel={Iterations},
    ytick={0,10,20,30,40,50,60},
    xtick={0.45,0.5,...,0.7},
    ymin=0,
    ymax=60,
    xlabel style={yshift=-0.cm},
    ylabel style={yshift=-0.35cm},
    xmajorgrids,
    ymajorgrids,
    legend style={at={(1,0.35)}, anchor=east, legend cell align=left, align=left, draw=white!15!black}
]
\addplot [color=TUDa-1a, line width=1pt,mark=o, mark options={}, mark size=2.0pt]
table[x=torque,y=iter,col sep=comma] {data/Pareto1st.txt};
\addlegendentry{First-order pareto}

\addplot [color=TUDa-8a, line width=1pt,mark=square, mark options={}, mark size=2.0pt]
table[x=torque,y=iter,col sep=comma] {data/Pareto2nd.txt};
\addlegendentry{Second-order pareto}

\end{axis}

\end{tikzpicture}
        \caption{Comparison of iterations for first-order vs.\ second-order solutions.}
        \label{fig:pareto:iteration}
    \end{subfigure}
    \caption{Pareto front generation for the PMSM case. First-order and second-order methods are compared in iterations. The second-order method clearly outperforms the first-order approach due to the good initial conditions.}
    \label{fig:pareto}
\end{figure*}
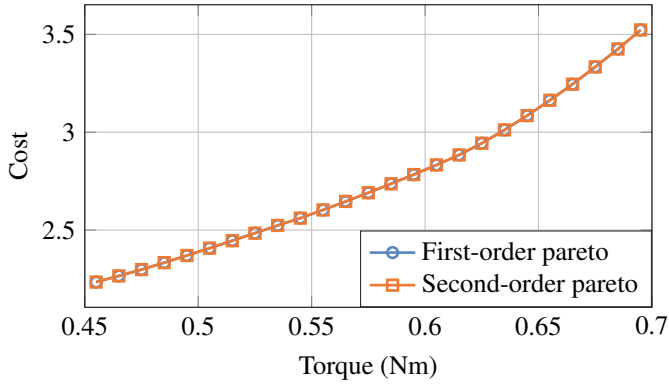
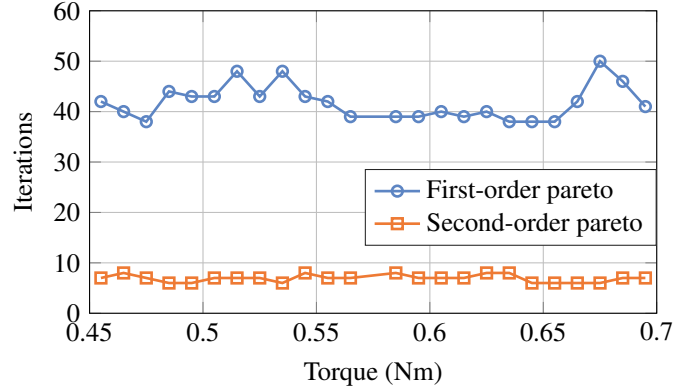

\subsection{Pareto front generation}
\label{MultiObjectiveOptimization}
Having established that the second-order Newton solver exhibits local superlinear convergence, improved search directions, and higher step quality for fixed scalar weighting of the objectives, we now take advantage of this property to trace the  Pareto front by a predictor-corrector continuation method based on the idea of \citep{Hillermeier_2001aa}. The resulting front answers the practical engineering question:  What is the cheapest rotor geometry that still achieves a prescribed average torque\;$\overline T$? 

The multi–objective continuation employed in this work is based on the $\varepsilon$-constraint method originally proposed by \citet{haimes2011multiobjective} in 1971. Specifically, we modify problem~\eqref{eq:FullOptimizationProblem} by introducing a user–selectable threshold~$\varepsilon_i$ for $i=1,2...m$ in the mean-torque constraint: 
\begin{equation}
    \label{eq:EpsProblem}
    \begin{aligned}
    \min_{\mathbf \design}\; F(\mathbf \design)
       &= m_{1}\,M(\mathbf \design)
          +m_{2}\,\widehat T(\mathbf \design)
          +m_{3}\,S(\mathbf \design), \\[2pt]
    \text{s.\,t.}\quad
    \mathbf K(\mathbf \state)\,\mathbf \state-\mathbf b &= \mathbf 0,\\
    \overline T(\mathbf \design) &\ge T_{\text{target}}-\varepsilon_i,\\
    \mathbf g(\mathbf \design) &\le \mathbf 0,
    \end{aligned}
\end{equation}
Note that alternative approaches also exist, such as adapting the weights from the weighted sum approach \citep{cesarano2024tracing}. We mention that the approach of \citet{cesarano2024tracing} can also recover nonconvex parts of the Pareto front, since it is a second-order method that computes all stationary points of the weighted-sum functional, including those that are not local minima. We nevertheless choose the $\varepsilon$-constraint method for the following reasons:
\begin{enumerate}
    \item \textbf{Physical Meaning.}  
          It is possible to specify the required average torque of the motor. When changing the weights, it is not known beforehand what the outcome will be. 
    \item \textbf{Optimum is attained at equality.}  
          In our model, using less magnet mass always lowers the torque.  At an optimum the limit is therefore reached when the constraint is sharp, $\overline T = T_{\text{target}}-\varepsilon_i$.
\end{enumerate}
Note that the proposed $\varepsilon$-constraint framework extends straightforwardly to higher-dimensional objective spaces by introducing additional thresholds for each quantity of interest.
The Pareto optimization starts by first solving \eqref{eq:EpsProblem} for an initial value~$\varepsilon_0$. For this, the optimization result from the previous section is used. Then, to obtain solutions for subsequent constraint values $\varepsilon_1,\varepsilon_2,\dots$, we simply use the previously converged solution as the initial guess for the new constraint. In practice, the step size \(\Delta\varepsilon\) can be chosen to ensure that neighboring constraint values remain close, so that each new subproblem is sufficiently close to the previously converged solution. This step size is dependent on the specific problem.  After each solver call, the new design $\design^{(k)}$ is stored and then the threshold is updated for the next run. The main advantage of this continuation approach is that the interior-point method starts from an almost-optimal geometry. As soon as the constraint is increased, it converges very quickly due to the locally superlinear behavior of the second-order Newton step.  In this way, a full sweep of \(\varepsilon\) values is obtained with minimal computational overhead compared to solving each instance from scratch. In our implementation, this continuation method corresponds to a zero-order predictor step, where the previously converged design is used directly as the initial guess for the next optimization run. While higher-order predictor–corrector schemes exploiting the tangent space of the Pareto manifold could yield improved initializations, the zero-order already achieves  sufficiently fast convergence.

This becomes also evident when applied to the PMSM case. For the Pareto front, we have relaxed the convergence criteria to a first order optimality of $10^{-3}$, because we are interested in the overall trend rather than a single optimum. This does not affect the key findings. In \cref{fig:pareto:cost}, we see two sets of Pareto curves plotted on the cost-versus-torque plane. Both methods result in the same curve. At lower torque levels, the relationship between cost and torque seems roughly linear. Once the torque demand exceeds a critical value, much more magnet mass is necessary, making the cost curve nonlinear. Such a knee point is usually interesting from a practical perspective to find a good tradeoff between cost and performance. The required iterations for both methods are shown in \cref{fig:pareto:iteration} and compared in \cref{tab:ParetoFront}. A similar shape of the Pareto curve was also found in \citep{Sato2022,COMSOL2021motoropt} who used a genetic algorithm.

\begin{table}[h]
    \caption{Comparison of average iterations, evaluations and computational time for the Pareto front generation of first- and second-order method.}
    \label{tab:ParetoFront}
    \centering
    \renewcommand{\arraystretch}{1.15}
    \begin{tabular}{|l||c|c|c|c|}
    \hline
    \textbf{Method} & \textbf{Iterations} & \textbf{Evaluations}  & \textbf{Runtime} \\ \hline\hline
    \textbf{First-order}& 41.8 & 77.4 &  \SI{6}{h} \SI{45}{min} \\ \hline
    \textbf{Second-order} & 6.9 & 8.3 &  \SI{4}{h} \SI{27}{min} \\ \hline
    \end{tabular}
\end{table}

\begin{figure*}[t]
    \centering
    \begin{subfigure}[b]{0.24\linewidth}
        \centering
        \includegraphics[width=\linewidth]{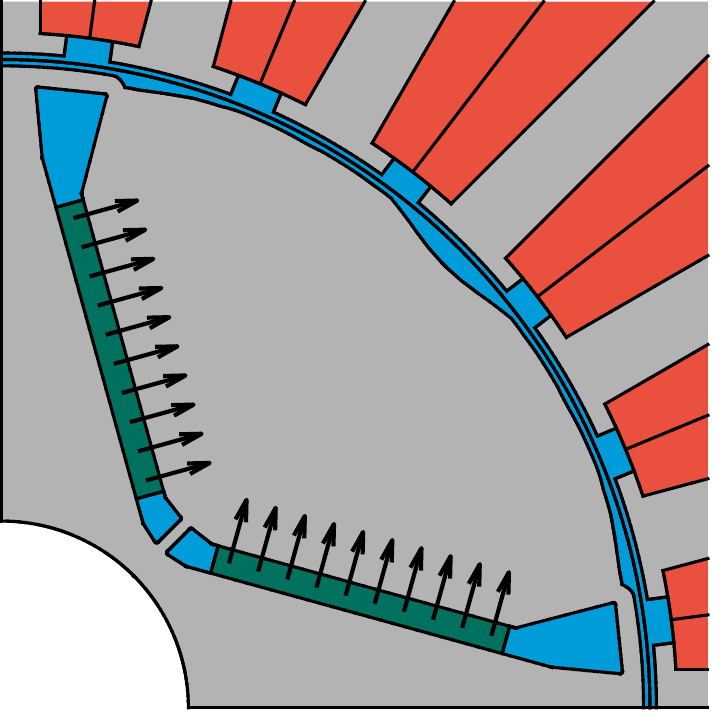}%
        \caption{Optimal geometry for $\overline{T}=\SI{0.455}{Nm}$.}
        \label{fig:ParetoGeometry1}
    \end{subfigure}%
    \hfill
    \begin{subfigure}[b]{0.24\linewidth}
        \centering
        \includegraphics[width=\linewidth]{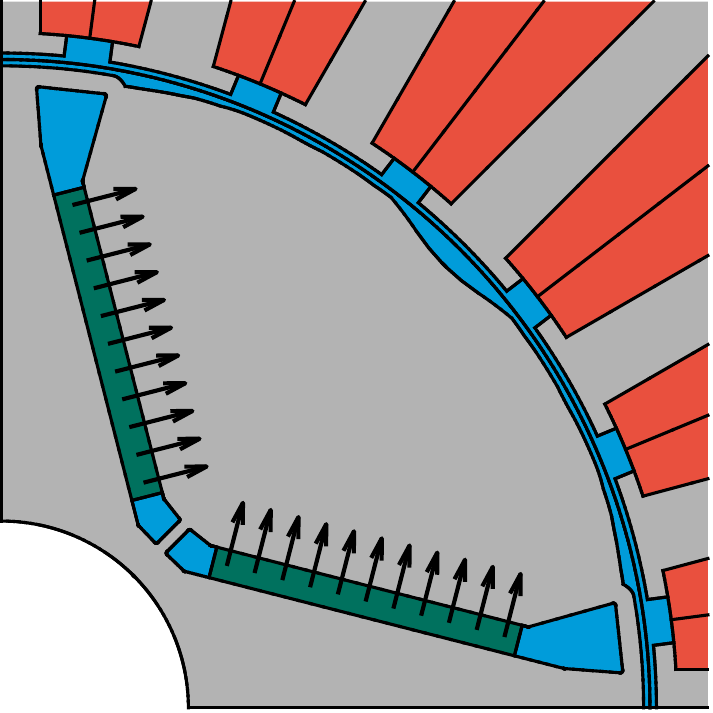}%
        \caption{Optimal geometry for $\overline{T}=\SI{0.515}{Nm}$.}
        \label{fig:ParetoGeometry2}
    \end{subfigure}%
    \hfill
    \begin{subfigure}[b]{0.24\linewidth}
        \centering
         \includegraphics[width=\linewidth]{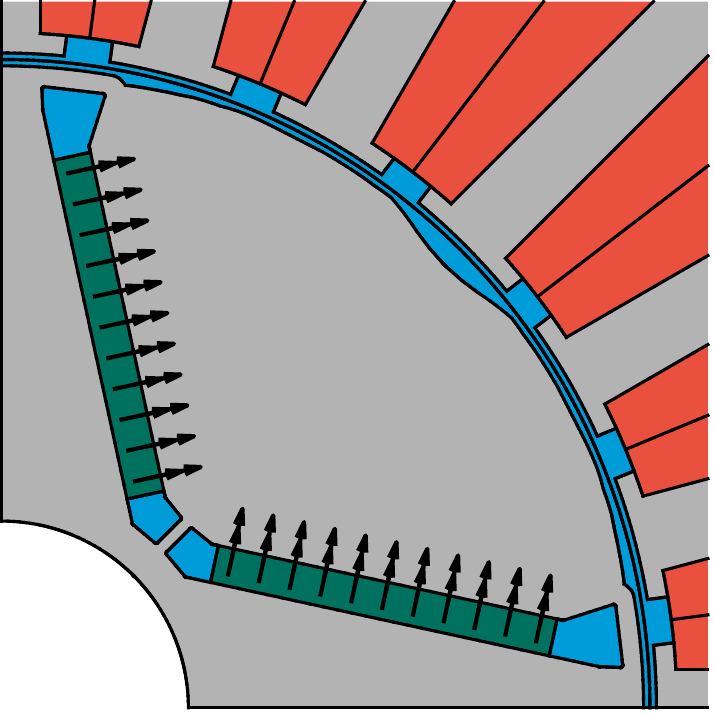}%
        \caption{Optimal geometry for $\overline{T}=\SI{0.625}{Nm}$.}
        \label{fig:ParetoGeometry3}
    \end{subfigure}%
    \hfill
    \begin{subfigure}[b]{0.24\linewidth}
        \centering
        \includegraphics[width=\linewidth]{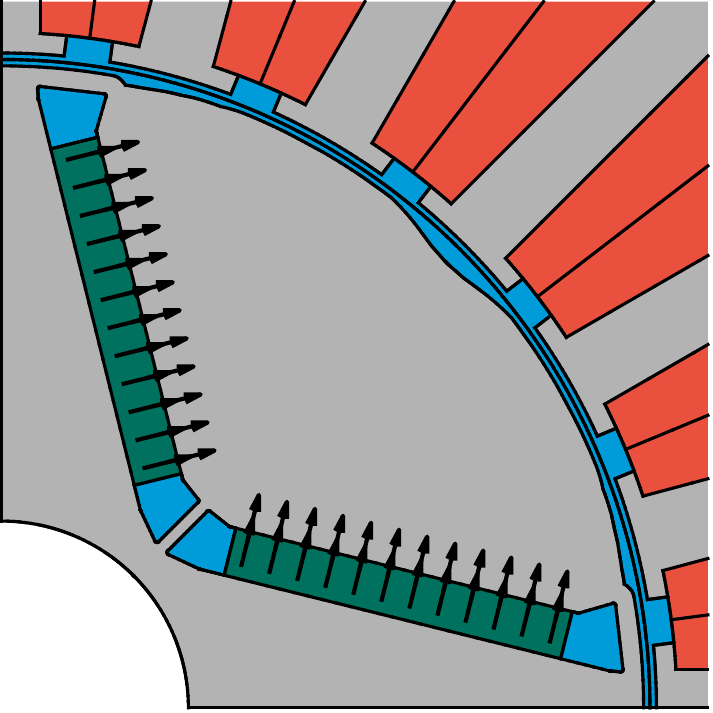}%
        \caption{Optimal geometry for $\overline{T}=\SI{0.695}{Nm}$.}
        \label{fig:ParetoGeometry4}
    \end{subfigure}

    \caption{Comparison of rotor geometries for different mean torque targets:
             (a)~optimized design at $\overline{T}=0.455$,
             (b)~optimized design at $\overline{T}=0.515$,
             (c)~optimized design at $\overline{T}=0.625$,
             (d)~optimized design at $\overline{T}=0.695$.}
    \label{fig:Pareto-Motors}
\end{figure*}

Finding all optima on the Pareto curve with the second-order method took \SI{4}{h}\,\SI{27}{min} with $6.9$ iterations on average. This is already considerably faster compared to the first-order method, which took \SI{6}{h}\,\SI{45}{min} with an average of $41.8$ iterations. 
The second-order continuation approach can exploit its locally superlinear convergence, has an improved search directions and a better step quality to move from one Pareto point to the next with minimal cost, because the new design is already close to an optimal solution from the previous run. In contrast, the first-order method does not benefit as much from previously calculated information. Here in particular, we see that second derivatives are valuable and can lower the computation time. 

A comparison of different motor geometries from the Pareto front is given in \cref{fig:Pareto-Motors}. The main difference between the designs is the size of the permanent magnets, because it majorly affects both cost and torque.

\section{Conclusions and Outlook}
\label{Conclusion}
In this paper, we derived second-order gradients for a design optimization problem constrained by the nonlinear magnetostatic PDE. The derivatives were applied for the design optimization of a permanent magnet synchronous machine, which is modeled with Isogeometric Analysis. This setting allows us to optimize the rotor parameters at the same time as the rotor surface. We show that the second-order method can significantly reduce the number of iterations compared to the first-order approach. In terms of computational runtime, the first-order method still outperforms the second-order method, at least  when optimizing from a standard initial design in the current framework. This changes when the optimization starts with a design that is already close to the optimum. Here, this is demonstrated by generating the cost-vs-torque Pareto front. In this case the locally superlinear convergence of the second-order method clearly outperforms the first-order approach both in iterations and runtime.

Prospective research areas include developing more efficient methods of calculating second derivatives, such as using Hessian-times-vector products or GPU-acceleration. Also extending the presented method to other physical aspects, for example to consider mechanical or thermal constraint already during the design process, would be valuable. Moreover, combining this framework with the existing robust optimization framework \citep{Komann2024,gangl2025robusttopologyoptimizationelectric} would be an interesting direction to identify robust Pareto points, since the current  designs require high accuracy, particularly in the air gap region.

\section*{Acknowledgement}

This work is supported by the joint DFG/FWF Collaborative Research Centre CREATOR (CRC -- TRR361/F90; DFG: Project-ID 492661287/TRR 361; FWF: Grant-DOI 10.55776/F90) at TU Darmstadt, TU Graz, RICAM and JKU Linz as well as the Graduate School CE within the Centre for Computational Engineering at TU Darmstadt. Moreover, the work of PG is partially supported by the State of Upper Austria.

\section*{Data Availability}
Link to Zenodo code repository will be inserted here upon acceptance.

\bibliography{reduced}

\appendix
\section{Jacobian matrix for Newton-Raphson}
\label{app:NewtonJacobian}
The derivative of the stiffness matrix $\stiffness_{ij}$ is given by
\begin{equation}
    \frac{\der\stiffness_{ij}}{\der \state_{k}} =\int _{\Omega }\underbrace{\frac{\partial \nu ( B)}{\partial B}}_{BH-curve}\frac{\der B}{\der \state_{k}} \nabla N_{i} \cdotp \nabla N_{j} \ \der\Omega 
\end{equation}
where $B=B( \state)$ and
\begin{align}
    \frac{\der B}{\der \state_{k}} & =\frac{\left( \nabla N_{l} \state_{l}\right) \cdotp \nabla N_{k}}{B}, & B= & \sqrt{\left( N_{l,x} \state_{l}\right)^{2} +\left( N_{l,y} \state_{l}\right)^{2}}. \label{eq:Bder}
\end{align}
In \eqref{eq:Bder} and also all following equations, the Einstein sum convention is used, usually over the index $l$. Indices after a comma indicate a derivative.

\section{First-order derivatives}
\label{app:first-order-derivatives}
The first order derivatives are obtained by differentiating the stiffness matrix \eqref{eq:methodology:Krt} and the right hand side vector \eqref{eq:methodology:brt} with respect to the control points $\Ckd$. Applying the product rule yields
\begin{equation}
    \frac{\der  \ \stiffness_{ij}}{\der  \ C_{kd}} =\stiffness_{ij,kd}^{( 1)} +\stiffness_{ij,kd}^{( 2)} +\stiffness_{ij,kd}^{( 3)} +\stiffness_{ij,kd}^{( 4)}
    \label{eq:app:firstOrderK}
\end{equation}
with the four components
\begin{align}
K_{ij,kd}^{( 1)} & =-\int_{\Omega} \nu(B)\,(\mathbf{D}_{kd} \nabla N_{i}) \cdotp \nabla N_{j}\ \der \Omega,\\
K_{ij,kd}^{( 2)} & =- \int_{\Omega} \nu(B)\,\nabla N_{i} \cdotp (\mathbf{D}_{kd} \nabla N_{j})\ \der \Omega,\\
K_{ij,kd}^{( 3)} & =\int_{\Omega} \nu(B)\,(\nabla N_{i} \cdot \nabla N_{j})\, \mathrm{Tr}(\mathbf{D}_{kd})\, \der \Omega,\\
K_{ij,kd}^{( 4)} &=
\int_{\Omega}
\underbrace{\frac{\partial \nu }{\partial B}}_{BH-curve}\frac{\der B}{\der C_{kd}} 
(\nabla N_i\cdot \nabla N_j)\, \der \Omega,
\end{align}
with 
\begin{equation}
    \mathbf{D}_{kd} =\begin{cases}
    \begin{pmatrix}
    \nabla G_{k} & \mathbf{0}
    \end{pmatrix} & \text{if}~~ d=x\\
    \begin{pmatrix}
    \mathbf{0} & \nabla G_{k}
    \end{pmatrix} & \text{if}~~d=y
    \end{cases}
\end{equation}
where $G$ comes from the physical mapping of the geometry functions in \eqref{eq:MappingFromParametricToPhysical}. All derivatives are derived on $\hat{\Omega}$ and then mapped back to $\Omega$. Hence integrals are written over $\Omega$ with $D_{kd}$ defined from $$G_k := \hat G_k\circ F^{-1}\qquad\text{on }\Omega$$
and
\begin{align}
    \frac{\der  B}{\der  C_{kd}} & =\frac{\left( \nabla N_{l} \state_{l}\right) \cdotp \left( -\mathbf{D}_{kd} \nabla N_{l} \state_{l}\right)}{\sqrt{\left( N_{l,x} \state_{l}\right)^{2} +\left( N_{l,y} \state_{l}\right)^{2}}}.
\end{align}

In the following we assume, that the coefficient $\nu$ appearing in $\rhs$ is constant
(or, more generally, independent of the design variables $C_{kd}$ in this term). The derivative of $\rhs$ with respect to the control points are 
\begin{equation}
    \frac{\der  \ \rhs_{i}}{\der  \ C_{kd}} =\rhs_{i,kd}^{( 1)} +\rhs_{i,kd}^{( 2)} \label{eq:app:firstOrderb}
\end{equation}
 with 
\begin{align}
\rhs_{i,kd}^{(1)}
&=
-\nu \int_{\Omega}
\big(\mathbf{D}_{kd}\,\nabla N_i\big)\cdotp \mathbf{B}_{r}^{\bot}\,\der\Omega,
\\[1mm]
\rhs_{i,kd}^{(2)}
&=
\nu \int_{\Omega}
\big(\nabla N_i\cdotp \mathbf{B}_{r}^{\bot}\big)\,
\mathrm{Tr}\big(\mathbf{D}_{kd}\big)\,\der\Omega .
\end{align}


\section{Second-order derivatives: $L_{\design\design}$}
\label{app:second-order-derivatives:Lxx}

The second derivatives of $\stiffness$ with respect to control points $\Ckd$ and $\Cle$ is obtained by differentiating \eqref{eq:app:firstOrderK} again with respect to $\Cle$. To improve readability we set $\nu=\nu(B)$. This yields the following contributions:
\begin{align}
K_{ij,kd,\ell e}^{( 1)( 1)} & =\int _{\Omega } \nu (\mathbf{D}_{\ell e}\mathbf{D}_{kd} \nabla N_{i}) \cdotp \nabla N_{j} \ \der \Omega \\
K_{ij,kd,\ell e}^{( 1)( 2)} & =\int _{\Omega} \nu (\mathbf{D}_{kd}\mathbf{D}_{\ell e} \nabla N_{i}) \cdotp \nabla N_{j} \ \der \Omega \\
K_{ij,kd,\ell e}^{( 1)( 3)} & =\int _{\Omega} \nu (\mathbf{D}_{kd} \nabla N_{i}) \cdotp (\mathbf{D}_{\ell e} \nabla N_{j}) \ \der \Omega\\
K_{ij,kd,\ell e}^{( 1)( 4)} & =\int _{\Omega} -\nu (\mathbf{D}_{kd} \nabla N_{i}) \cdotp \nabla N_{j} \  \mathrm{Tr}(\mathbf{D}_{\ell e}) \ \der \Omega \\
K_{ij,kd,\ell e}^{( 1)( 5)} & =\int _{\Omega } -\frac{\partial \nu }{\partial B}\frac{\der B}{\der C_{\ell e}}(\mathbf{D}_{kd} \nabla N_{i}) \cdotp \nabla N_{j} \  \der \Omega 
\end{align}

\begin{align}
K_{ij,kd,\ell e}^{( 2)( 1)} & =\int _{\Omega} \nu (\mathbf{D}_{\ell e} \nabla N_{i}) \cdotp (\mathbf{D}_{kd} \nabla N_{j}) \ \der \Omega\\
K_{ij,kd,\ell e}^{( 2)( 2)} & =\int _{\Omega } \nu \nabla N_{i} \cdotp (\mathbf{D}_{\ell e}\mathbf{D}_{kd} \nabla N_{j}) \ \der \Omega\\
K_{ij,kd,\ell e}^{( 2)( 3)} & =\int _{\Omega } \nu \nabla N_{i} \cdotp (\mathbf{D}_{kd}\mathbf{D}_{\ell e} \nabla N_{j}) \ \der \Omega\\
K_{ij,kd,\ell e}^{( 2)( 4)} & =\int _{\Omega } -\nu \nabla N_{i} \cdotp (\mathbf{D}_{kd} \nabla N_{j}) \ \mathrm{Tr}(\mathbf{D}_{\ell e}) \ \der \Omega \\
K_{ij,kd,\ell e}^{( 2)( 5)} & =\int _{\Omega} -\frac{\partial \nu }{\partial B}\frac{\der B}{\der C_{\ell e}} \nabla N_{i} \cdotp (\mathbf{D}_{kd} \nabla N_{j}) \ \der \Omega 
\end{align}

\begin{align}
K_{ij,kd,\ell e}^{( 3)( 1)} & =\int _{\Omega} -\nu (\mathbf{D}_{\ell e} \nabla N_{i}) \cdotp \nabla N_{j} \ \ \mathrm{Tr}(\mathbf{D}_{kd})  \der \Omega\\
K_{ij,kd,\ell e}^{( 3)( 2)} & =\int _{\Omega } -\nu \nabla N_{i} \cdotp (\mathbf{D}_{\ell e} \nabla N_{j}) \ \mathrm{Tr}(\mathbf{D}_{kd}) \ \der \Omega\\
K_{ij,kd,\ell e}^{( 3)( 3)} & =\int _{\Omega} \nu \nabla N_{i} \cdotp \nabla N_{j} \  \mathrm{Tr}(\mathbf{D}_{\ell e}) \ \mathrm{Tr}(\mathbf{D}_{kd}) \ \der \Omega\\
K_{ij,kd,\ell e}^{( 3)( 4)} & =\int _{\Omega} -\nu \nabla N_{i} \cdotp \nabla N_{j} \  \mathrm{Tr}(\mathbf{D}_{\ell e}\mathbf{D}_{kd}) \ \der \Omega \\
K_{ij,kd,\ell e}^{( 3)( 5)} & =\int _{\Omega}\frac{\partial \nu }{\partial B}\frac{\der B}{\der C_{\ell e}} \nabla N_{i} \cdotp \nabla N_{j} \ \mathrm{Tr}(\mathbf{D}_{kd}) \ \der \Omega 
\end{align}

\begin{align}
K_{ij,kd,\ell e}^{( 4)( 1)} & =\int _{\Omega} -\frac{\partial \nu }{\partial B}\frac{\der B}{\der C_{kd}}(\mathbf{D}_{\ell e} \nabla N_{i}) \cdotp \nabla N_{j} \  \der \Omega \\
K_{ij,kd,\ell e}^{( 4)( 2)} & =\int _{\Omega} -\frac{\partial \nu }{\partial B}\frac{\der B}{\der C_{kd}} \nabla N_{i} \cdotp (\mathbf{D}_{\ell e} \nabla N_{j})  \  \der \Omega\\
K_{ij,kd,\ell e}^{( 4)( 3)} & =\int _{\Omega}\frac{\partial \nu }{\partial B}\frac{\der B}{\der C_{kd}} \nabla N_{i} \cdotp \nabla N_{j} \ \mathrm{Tr}(\mathbf{D}_{\ell e}) \  \der \Omega\\
K_{ij,kd,\ell e}^{( 4)( 4)} & =\int _{\Omega}\frac{\partial \nu }{\partial B}\frac{\der B}{\der C_{kd} \der C_{\ell e}} \nabla N_{i} \cdotp \nabla N_{j}  \der \Omega\\
K_{ij,kd,\ell e}^{( 4)( 5)} & =\int _{\Omega}\frac{\partial ^{2} \nu }{\partial B^{2}}\frac{\der B}{\der C_{\ell e}}\frac{\der B}{\der C_{kd}} \nabla N_{i} \cdotp \nabla N_{j} \  \der \Omega
\end{align}
with

\begin{align}
\frac{\der B}{\der C_{kd} \der C_{\ell e}} & =\frac{\left(\mathbf{D}_{\ell e} \nabla N_{l} \state_{l}\right) \cdotp \left(\mathbf{D}_{kd} \nabla N_{l} \state_{l}\right)}{B} \nonumber\\
 & +\frac{\left( \nabla N_{l} \state_{l}\right) \cdotp \left(\mathbf{D}_{\ell e}\mathbf{D}_{kd} \nabla N_{l} \state_{l}\right)}{B}\\
 & +\frac{\left( \nabla N_{l} \state_{l}\right) \cdotp \left(\mathbf{D}_{kd}\mathbf{D}_{\ell e} \nabla N_{l} \state_{l}\right)}{B} \nonumber\\
 & -\frac{1}{B}\frac{\der B}{\der C_{kd}}\frac{\der B}{\der C_{\ell e}} \nonumber
\end{align}

The same is achieved for the right hand side. Note, that $\nu$ is here constant. By differentiation of \eqref{eq:app:firstOrderb}, we obtain the following contributions:

\begin{align}
    b_{i,kd,\ell e}^{( 1)( 1)} = & \nu \int _{\Omega}(\mathbf{D}_{\ell e}\mathbf{D}_{kd} \nabla N_{i}) \cdotp \mathbf{B}_{\mathrm{r}}^{\bot }  \ \der \Omega\\
    b_{i,kd,\ell e}^{( 1)( 2)} = & \nu \int _{\Omega}(\mathbf{D}_{kd}\mathbf{D}_{\ell e} \nabla N_{i}) \cdotp \mathbf{B}_{\mathrm{r}}^{\bot } \  \der \Omega\\
    b_{i,kd,\ell e}^{( 1)( 3)} = & \nu \int _{\Omega} -(\mathbf{D}_{kd} \nabla N_{i}) \cdotp \mathbf{B}_{\mathrm{r}}^{\bot }  \mathrm{Tr}(\mathbf{D}_{\ell e}) \ \der \Omega
\end{align}

\begin{align}
b_{i,kd,\ell e}^{( 2)( 1)} = & \nu \int _{\Omega } -(\mathbf{D}_{\ell e} \nabla N_{i}) \cdotp \mathbf{B}_{\mathrm{r}}^{\bot }  \mathrm{Tr}(\mathbf{D}_{kd}) \ \der \Omega\\
b_{i,kd,\ell e}^{( 2)( 2)} = & \nu \int _{\Omega } \nabla N_{i} \cdotp \mathbf{B}_{\mathrm{r}}^{\bot } \ \mathrm{Tr}(\mathbf{D}_{\ell e}) \ \mathrm{Tr}(\mathbf{D}_{kd}) \  \der \Omega\\
b_{i,kd,\ell e}^{( 2)( 3)} = & \nu \int _{\Omega } -\nabla N_{i} \cdotp \mathbf{B}_{\mathrm{r}}^{\bot } \ \mathrm{Tr}(\mathbf{D}_{\ell e}\mathbf{D}_{kd}) \ \der \Omega 
\end{align}

The components for $L_{xx}$ are then obtained by multiplying these derivatives with adjoint and state vector.

\section{Second-order derivatives: $L_{\design\state}$}
\label{app:second-order-derivatives:Lxy}
The mixed order derivatives are obtained by differentiation with respect to $\design$ and $\state$. Only the nonlinear part of the stiffness matrix yields new contributions in this case:
\begin{align}
\frac{\der }{\der C_{kd}}\left(\frac{dK}{dy} \cdot y\right)_{ij} =& -\int _{\hat{\Omega }}\frac{\partial \nu }{\partial B}\frac{\der B}{\der \state_{j}}\left(\mathbf{D}_{kd} \nabla N_{i}\right) \cdot \nabla N_{l} \state_{l}  \ \der \Omega \\
 & -\int _{\hat{\Omega }}\frac{\partial \nu }{\partial B}\frac{\der B}{\der \state_{j}} \nabla N_{i} \cdot \left(\mathbf{D}_{kd} \nabla N_{l} \state_{l} \right) \  \der \Omega \\
 & +\int _{\hat{\Omega }}\frac{\partial \nu }{\partial B}\frac{\der B}{\der \state_{j}} \nabla N_{i} \cdot \nabla N_{l} \state_{l}\mathrm{Tr}(\mathbf{D}_{kd}) \ \der \Omega \\
 & +\int _{\hat{\Omega }}\frac{\partial ^{2} \nu }{\partial B^{2}}\frac{\der B}{\der C_{kd}}\frac{\der B}{\der \state_{j}} \nabla N_{i} \cdot \nabla N_{l} \state_{l}  \  \der \Omega\\
 & +\int _{\hat{\Omega }}\frac{\partial \nu }{\partial B}\frac{\der ^{2} B}{\der \state_{j}\der C_{kd}} \nabla N_{i} \cdot \nabla N_{k} \state_{l} \   \der \Omega 
\end{align}

with the mixed term derivative term
\begin{align}
\frac{\der ^{2} B}{\der \state_{j}\der C_{kd}} & =-\frac{ (\mathbf{D}_{kd} \nabla N_{l} \state_{l}) \cdot \nabla N_{j}}{B} \nonumber\\
 & -\frac{ ( \nabla N_{l} \state_{l}) \cdot \mathbf{D}_{kd} \nabla N_{j}}{B}\\
 & -\frac{ ( \nabla N_{l} \state_{l}) \cdot \nabla N_{j}}{B^{2}}\frac{\der B}{\der C_{kd}}. \nonumber
\end{align}

\section{Second-order derivatives: $L_{\state\state}$}
\label{app:second-order-derivatives:Lyy}

The component $L_{\state\state}$ has different components. First, the torque \eqref{eq:Torque} is differentiated twice, resulting in
\begin{align}
    L_{yy}^{( 1)} & =-l_z\begin{pmatrix}
    0 & 0 & 0\\
    0 & 0 & G_{\mathrm{st}} R_{\alpha }^{'}\\
    0 & \left( G_{\mathrm{st}} R_{\alpha }^{'}\right)^{\top } & 0
    \end{pmatrix}
\end{align}

In addition, again the nonlinear part of the stiffness matrix does not vanish. The contribution is 
\begin{multline}
    \begin{aligned}
    L_{yy}^{(2)} =  -\int_{\Omega}z_{l} & \Biggl(\left(\frac{\partial ^{2} \nu }{\partial B^{2}}\frac{\der B}{\der \state_{i}}\frac{\der B}{\der \state_{j}} +\frac{\partial \nu }{\partial B}\frac{\der ^{2} B}{\der \state_{i}\der \state_{j}}\right) \nabla N_{l} \cdot \nabla N_{m} \state_{m}\\
     & +\frac{\partial \nu }{\partial B}\frac{\der B}{\der \state_{i}} \nabla N_{l} \cdot \nabla N_{j}\\
     & +\frac{\partial \nu }{\partial B}\frac{\der B}{\der \state_{j}} \nabla N_{l} \cdot \nabla N_{i}\Biggl) \ \der \Omega
    \end{aligned}
\end{multline}
with 
\begin{align}
\frac{\der ^{2} B}{\der \state_{i}\der \state_{j}} & =\frac{\nabla N_{i} \cdot \nabla N_{j}}{B} -\frac{ \nabla N_{l} \state_{l} \cdot \nabla N_{j}}{B^{2}}\frac{\der B}{\der \state_{i}}.
\end{align}

\end{document}